\documentclass[preprint,12pt,times]{elsarticle}

\usepackage{amsmath,amssymb,amsfonts}
\usepackage[justification=centering]{caption}
\usepackage{xpatch}
\usepackage{subfigure}
 \usepackage{pdfsync}
\usepackage{braket,amsfonts}
\usepackage{bm,array}
\usepackage{placeins,graphicx,diagbox }
\usepackage{etex,amsmath,amssymb}  
\usepackage{pgfplots,stfloats}
\usepackage{geometry}
\usepackage{amsmath} 
\usepackage{amssymb}  
\usepackage{color}
\usepackage{graphicx}
\usepackage{caption}
\usepackage{graphicx}
\usepackage{float}
\usepackage{bm}
\usepackage{color}
\usepackage{framed}
\usepackage{color}
\usepackage{listings}
\usepackage{multicol}
\usepackage{tikz}

\newtheorem{algorithm}{Algorithm}
\usepackage{algorithm}
\usepackage{algorithmic}
\usepackage{bbm,extarrows}
\usepackage{multirow,array}

\makeatletter
\def\changeBibColor#1{%
	\in@{#1}{}
	\ifin@\color{blue}\else\normalcolor\fi
}

\xpatchcmd\@bibitem
{\item}
{\changeBibColor{#1}\item}
{}{\fail}

\xpatchcmd\@lbibitem
{\item}
{\changeBibColor{#2}\item}
{}{\fail}
\makeatother
\usepackage{hyperref}
\hypersetup{
	colorlinks=true,
	linkcolor=blue,
	filecolor=blue,
	urlcolor=blue,
	citecolor=cyan,
}
\def\argmin{\mathop{\rm argmin}}
\newtheorem{theorem}{Theorem}
\newtheorem{lemma}{Lemma}

\newtheorem{assumption}{Assumption}
\newtheorem{proposition}{Proposition}

\biboptions{numbers,sort&compress}

\begin{document}

\begin{frontmatter}



\title{Distributed Projection-free Algorithm for Constrained Aggregative Optimization}
\tnotetext[t1]{This research is supported by  Shanghai Municipal Science and Technology Major Project No. 2021SHZDZX0100,  Shanghai Municipal Commission of Science and Technology Project No. 19511132101,
	National Natural Science Foundation of China under Grant 62003239, Shanghai Sailing Program No. 20YF1453000,
	and the Fundamental Research Funds for the Central Universities, China.}
\author[rvt,gf]{Tongyu Wang}
\ead{wang\_tongyu@tongji.edu.cn}
\author[rvt,ai,gf]{Peng Yi }
\ead{yipeng@tongji.edu.cn}

\address[rvt]{Department of Control Science and Engineering, Tongji University,  Shanghai,   201804, China;}
\address[ai]{The Shanghai Research
	Institute for Intelligent Autonomous Systems,Shanghai, 201804, China.}
\address[gf]{Shanghai Institute of Intelligent Science and Technology, Tongji University, Shanghai 200092, China}


\begin{abstract}
In this paper, we focus on solving a distributed convex aggregative optimization problem in a network, where each agent has its own cost function which depends  not only on its own decision variables but also on the aggregated function of all agents' decision variables. The decision variable is constrained within a feasible set. In order to minimize the sum of the cost functions when each agent only knows its local cost function, we propose a distributed Frank-Wolfe algorithm based on gradient tracking for the aggregative optimization problem where each node maintains two estimates, namely an estimate of the sum of agents' decision variable and an estimate of the gradient of global function. The algorithm is projection-free, but only involves solving a linear optimization to get a search direction at each step.
We show the convergence of the proposed algorithm for convex and smooth objective functions over a time-varying network. Finally, we demonstrate the convergence and computational efficiency of the proposed algorithm via numerical simulations.
\end{abstract}



\begin{keyword}
Distributed algorithm,
Aggregative optimization,
Gradient projection-free,
Time-varying graph


\end{keyword}

\end{frontmatter}


\section{Introduction}\label{int}
In recent years, distributed optimization algorithms and their applications have received extensive attention in  decision-making problems for multi-agent networks\cite{nedic2018distributed,notarstefano2019distributed,yang2019survey}, with applications in smart grids \cite{Nguyen19}, resource allocation \cite{DAI2022}, and robot formations\cite{BHOWMICK2022}. In the framework of multi-agent networks, the agents in the network have a local interactive  pattern, each of them can only access its own information and that of  its neighboring agents, and the goal of the agents is to optimize the sum of the local objective functions in a cooperative manner.

In general, distributed optimization algorithms can be divided into unconstrained optimization \cite{nedic2009distributed} \cite{chen2020distributed}, and constrained optimization \cite{nedic2018distributed} \cite{yi2015distributed},  from the viewpoint of  with or without constraints. For unconstrained optimization, various algorithms like consensus  subgradient algorithm \cite{nedic2009distributed}, dual averaging \cite{duchi2011dual}, EXTRA \cite{shi2015extra} and gradient tracking algorithm \cite{pu2021distributed} are studied.
For constrained optimization problems, methods based on projection dynamics and primal-dual dynamics have been  studied. For example, \cite{nedic2014distributed} studied the distributed projected subgradient algorithm. \cite{lei2016primal} studied a projected primal-dual algorithm for constrained optimization, and  \cite{Liang2020} developed a distributed  dual average push-sum algorithm with dual decomposition, while a constrained subproblem is solved at each step. \cite{YUAN2018} investigated a distributed mirror-descent optimization algorithm based on the Bregman divergence and achieved an $O(ln(T)/T)$ rate of convergence. \cite{Yuan2016} proposed a consensus-based distributed regularized primal-dual gradient method. Compared to  methods that require projection of estimates onto the constraint set at each iteration, the algorithm in \cite{Yuan2016} only required one projection at the last iteration.

However, the projection-based distributed algorithm implies that the agent needs to solve  a quadratic optimization problem at each iteration, to find the closest point within constraint set. When the constraints have complex structures (e.g., polyhedra), the computational cost of solving quadratic subproblem can prevent the agent from using projection-based dynamics, especially for high-dimensional optimization problems.  Fortunately, the well-known Frank-Wolfe(FW) \cite{frank1956algorithm} algorithm provides us with a way to derive an effective searching direction while maintaining decision feasibility. Each step of FW algorithm only needs to solve a constrained linear programming problem, which could have a closed form for specific problems or have effective
solvers.
Recently, FW methods have received renewed research attention due to its projection-free property and advantages for large-scale problems for online learning \cite{hazan2020faster}, machine learning \cite{chen2018projection} and traffic assignment \cite{chen2001effects}.
Briefly speaking, the FW method uses a linearized function to approximate the objective function and derives feasible descent directions by solving a linear objective optimization, by $\boldsymbol{s}:=\underset{\boldsymbol{s} \in \mathcal{D}}{\arg \min }\left\langle\boldsymbol{s}, \nabla f\left(\boldsymbol{x}^{(k)}\right)\right\rangle$, $\boldsymbol{x}^{(k+1)}:=(1-\gamma) \boldsymbol{x}^{(k)}+\gamma \boldsymbol{s}$. There have been massive developments and applications for the FW method afterwards. For example,  the primal-dual convergence rate has been analyzed in \cite{jaggi2013revisiting}, and \cite{lacoste2016convergence} have analyzed its convergence for non-convex problems. In addition, \cite{Lafond2016} developed communication efficient algorithms using the stochastic Frank-Wolfe (sFW) algorithm when there exist noises on the gradient computation. FW algorithms are also utilized and extended to distributed optimization. For example, \cite{wai2017decentralized} studied the decentralized optimization problem by treating the algorithm as an inexact FW algorithm with consensus errors.   By combing gradient tracking and FW algorithm, \cite{Guanpu2021} investigated the continuous-time algorithm based on FW dynamics, and \cite{jiang2022distributed} studied FW algorithm for constrained stochastic optimization.

In the distributed optimization problems mentioned above, the agents decide the same variable while need to reach consensus on the optimal decision. However, in other  scenarios, each agent only decides its own variable but the objective function of each agent is related to other agents' decisions through an aggregated variable. For example, in multi-agents formation control problem \cite{cao2020distributed},  a group of networked agents wish to achieve a geometric pattern while surrounding a target, which can be regarded as a goal tracking problem. The dynamical tracking problem can be modelled as an online distributed optimization.  In this case, each agent decides its location, while the objective function also depends on the centroid of all agents. Similar scenarios exist in resource allocation \cite{yi2016initialization}, smart grids \cite{Longe2017}, social networks \cite{PANTOJA2019209} as well. The above optimization problem is called distributed aggregative optimization in \cite{li2021distributed}. The aggregative optimization is also related to the well studied aggregative games \cite{koshal2016distributed,liang2017distributed,grammatico2017dynamic}, but can be treated as a cooperative formulation in contrast to the noncooperative setting for multiagent decision problem. Regarding the study of aggregative optimization, \cite{li2021distributed} considered a static unconstrained framework, \cite{carnevale2021distributed} considered an online constrained framework, and \cite{chen2022distributed} considered a quantitative problem. To the best of our knowledge, no work has been done to improve the computational bottlenecks encountered by algorithms when dealing with complicated constraints that prohibiting projection algorithms.

Therefore, our interest is to design projection-free methods to solve distributed aggregative optimization problems with constraints. Motivated both by \cite{li2021distributed} and \cite{jaggi2013revisiting}, we propose a FW based approach to solve the aggregative optimization in a distributed manner. The main contributions are as follows.
\begin{itemize}
\item Firstly, a novel distributed projection-free algorithm based on FW method  with gradient tracking is designed to solve the aggregative optimization problem. Each agent's local cost function depends both on its own variable and on aggregated variable, while the global information are not known by any single agent. The proposed algorithm uses the dynamical averaging tracking approach to estimate the global aggregation variable and the corresponding gradient term.
\item Secondly, we prove that the algorithm converges to the optimal solution when the objective function is  convex. Compared with the projected dynamics in \cite{carnevale2021distributed}, the proposed algorithm is able to solve the aggregative optimization problem  over  time-varying communication graphs.
\item Finally, we  demonstrate the efficiency of the proposed algorithm with numerical studies.
\end{itemize}
The  rest of paper is organized as follows. Section 2 introduces  notations, preliminary knowledge of graph theory, and illustrates the distributed aggregative optimization problem. Section 3 provides the proposed distributed algorithm and main convergence result. Section 4 provides the proof of convergence. Then, a numerical experiment is given in Section 5 and Section 6 concludes the paper.

{\bf Notations:} When referring to a vector $x$, it is assumed to be
a column vector while $x^{\top}$ denotes its transpose. $\langle  x,y\rangle=x^{\top}y$ denotes the inner product of vectors $x,y.$  $\|x\|$
{denotes} the Euclidean vector norm, i.e., $\|x\|=\sqrt{x^{\top}x}$. Let $\otimes$ be the Kronecker product. Denote by $\mathbf{1}_{N}$ and $\mathbf{0}_{N}$ the column vectors.
A nonnegative square   matrix $A $   is  called doubly stochastic if
$A\mathbf{1} =\mathbf{1}$ and  $\mathbf{1}^{\top} A =\mathbf{1}^{\top}$, where
$\mathbf{1}$ denotes the vector with each entry being $1$.
$\mathbf{I}_N \in \mathbb{R}^{N\times N}$ denotes the  identity matrix.
Let $\mathcal{G}=\{ \mathcal{N}, \mathcal{E}\}$ be  a directed graph   with
$\mathcal{N}=\{1,\cdots,N\} $  denoting the set of  players and  $\mathcal{E}$
denoting  the set of  directed edges between players, where
$(j,i)\in\mathcal{E }$ if  player  $i$ can  obtain  information from
player  $j$.
The graph  $\mathcal{G}$ is  called  strongly connected if  for
any $  i,j\in \mathcal{N}$ there exists a directed path from  $i$ to
$j$, i.e., there exists  a  sequence  of edges   $ (i,i_1),(i_1,i_2),\cdots,(i_{p-1},j)$ in the digraph  with  distinct  nodes $ i_m \in \mathcal{N},~\forall m: 1 \leq m \leq p-1$.
A differentiable function $f: \mathbb{R}^{n} \rightarrow \mathbb{R}$ is $\mu$-strongly convex if for all $\boldsymbol{\theta}, \boldsymbol{\theta}^{\prime} \in \mathbb{R}^{n}$,$
f(\boldsymbol{\theta})-f\left(\boldsymbol{\theta}^{\prime}\right) \leq\left\langle\nabla f(\boldsymbol{\theta}), \boldsymbol{\theta}-\boldsymbol{\theta}^{\prime}\right\rangle-\frac{\mu}{2}\left\|\boldsymbol{\theta}-\boldsymbol{\theta}^{\prime}\right\|_{2}^{2}
$. moreover, $f$ is convex if the above is satisfied with $\mu=0$.

\section{Problem  Formulation }\label{sec:formulation}

In this section, we  formulate the  aggregative optimization problem over networks and introduce some basic assumptions.

\subsection{Problem Statement}

We define the  aggregative optimization problem with $N$ agents:
\begin{equation}\label{Nopt_agg}
\begin{split}
	\min_{\boldsymbol{x} \in X \triangleq  \prod_{i=1}^{N} {X}_{i}}& f(\boldsymbol{x}) = \sum_{i=1}^{N}f_{i}(x_{i},\delta(\boldsymbol{x}))  {\rm~with~}\delta(\boldsymbol{x}) = \frac{1}{N}\sum_{i=1}^{N}\phi_{i}(x_{i}),
\end{split}
\end{equation}
where $\boldsymbol{x} \triangleq \operatorname{col}\big((x_{i})_{i\in \mathcal{N}}\big)$ is the global strategy variable with $x_{i}\in\mathbb{R}^{n_{i}}$ being the decision variable of agent $i$. The function $f_{i}:\mathbb{R}^{n}\rightarrow\mathbb{R}$ is the local objective function of agent $i$ with $n\triangleq\sum_{i=1}^{N}n_{i}$,  $X_i \subset \mathbb{R}^{n_i}$ denotes the local constraint set of agent $i,$  and $\mathcal{N} =\{1,\dots,N\}$ denotes the set of  agents.    Moreover,  $\delta(\boldsymbol{x})$ is an aggregate information of all agents' variables, where $\phi_{i}:\mathbb{R}^{n_i}\rightarrow\mathbb{R}^{d}$ is only available to agent $i$.
The goal is to design a distributed algorithm to cooperatively seek an optimal decision variable to the problem 
(\ref{Nopt_agg}).

 The gradient of $f(\boldsymbol{x})$ is defined by
\[\nabla f(\boldsymbol{x}) \triangleq \nabla_{1} f(\boldsymbol{x}, \boldsymbol{\delta}(\boldsymbol{x}))+\nabla \phi(\boldsymbol{x}) \mathbf{1}_{N} \otimes \frac{1}{N} \sum_{i=1}^{N} \nabla_{\delta} f_{i}\left(x_{i},\delta(\boldsymbol{x})\right),
\] where $\nabla_{1} f(\boldsymbol{x},\boldsymbol{\delta}(\boldsymbol{x})): =\operatorname{col} (\nabla_{x_i}f_{i}(x_{i},\mathbf{z})_{|\mathbf{z}=\boldsymbol{\delta}(\boldsymbol{x})})_{i\in\mathcal{N}} $, and
\[\nabla \phi\left(\boldsymbol{x}\right):=\left[\begin{array}{lll}
	\nabla \phi_{1}\left(x_{1}\right) & & \\
	& \ddots & \\
	& & \nabla \phi_{N}\left(x_{N}\right)
\end{array}\right] .\]

For ease of notation, we also specify the cost function as $f_{i}(x_{i},\delta(\boldsymbol{x})) = g_{i}(x_{i},z_{i})_{|z_{i}=\delta(\boldsymbol{x})}$ with a function $g_{i}:\mathbb{R}^{n_{i}+n} \rightarrow \mathbb{R}, i\in \mathcal{N}$. To move forward, define $g(\boldsymbol{x},\boldsymbol{z}) = \sum_{i=1}^{N}g_{i}(x_{i},z_{i}):\mathbb{R}^{n+Nd}\rightarrow \mathbb{R}$ for any $\boldsymbol{x}$ and $\boldsymbol{z} \triangleq \operatorname{col}((z_{i})_{i\in \mathcal{N}} )$.  The gradient of $g(\boldsymbol{x},\boldsymbol{z})$ is defined by $\nabla_{1}g(\boldsymbol{x},\boldsymbol{z}):=\operatorname{col}(\nabla_{x_i}g_{i}(x_{i},z_{i})_{i\in\mathcal{N}})$, $\nabla_{2}g(\boldsymbol{x},\boldsymbol{z}):=\operatorname{col}(\nabla_{z_i}g_{i}(x_{i},z_{i})_{i\in\mathcal{N}})$.

Next, we impose some  assumptions on the formulated problem.
We require the agent-specific problem to be convex and  continuously differentiable.

\begin{assumption}~\label{ass-payoff}
	(a) For each agent $i\in \mathcal{N},$ the strategy set  $X_i$ is closed, convex and compact.
 In addition, the diameter of $X = \prod_{i=1}^{N} {X}_{i} $ is defined as $\bar{\rho}:=\max_{\theta,\theta^{\prime}\in X}\|\theta-\theta^{\prime}\|^{2}_{2}$;
	\\ (b) the global objective function $f$ is convex and differentiable in $\boldsymbol{x} \in X$, and
 its gradient function is  $L$-smooth in $\boldsymbol{x} \in X$, i.e.,
	\[ f(\mathbf{y})-f(\mathbf{x}) \geq (\mathbf{y}-\mathbf{x})^T \nabla f(\mathbf{x}),~\| \nabla f(\mathbf{x}) -\nabla f(\mathbf{y}) \|\leq L\|\mathbf{x}-\mathbf{y}\|, \quad \forall \mathbf{x},\mathbf{y} \in X.\]
\end{assumption}

\begin{assumption}~\label{ass-gradient}
	(a) $\nabla_{1} g(\boldsymbol{x}, \boldsymbol{z})$ is $l_{1}$-Lipschitz continuous in $\boldsymbol{z} \in \mathbb{R}^{Nd}$ for any $\boldsymbol{x} \in X,$   i.e.,
\[ \|\nabla_{1} g(\boldsymbol{x}, \boldsymbol{z}_1)-\nabla_{1} g(\boldsymbol{x}, \boldsymbol{z}_1)\| \leq l_1 \|\boldsymbol{z}_1-\boldsymbol{z}_2\|,~\forall \boldsymbol{z}_1, \boldsymbol{z}_2 \in \mathbb{R}^{Nd}.\]
	(b) $\nabla_{2} g(\boldsymbol{x}, \boldsymbol{z})$ is $l_{2}$-Lipschitz continuous in $(\boldsymbol{x}, \boldsymbol{z}) \in  X \times \mathbb{R}^{Nd}$, i.e.
	\begin{align}
		\|\nabla_{2} g(\boldsymbol{x}_1, \boldsymbol{z}_1)-\nabla_{2} g(\boldsymbol{x}_2, \boldsymbol{z}_2)\| \leq l_2 \|\boldsymbol{x}_1-\boldsymbol{x}_1\| +l_2 \|\boldsymbol{z}_1-\boldsymbol{z}_2\|, \forall \boldsymbol{x}_1, \boldsymbol{x}_2 \in X, \boldsymbol{z}_1, \boldsymbol{z}_2 \in \mathbb{R}^{Nd}.
	\end{align}
	(c) For each agent $i\in \mathcal{N}$, $\phi_{i}(x_{i})$ is differentiable in $x_i \in X_i$ and $\nabla\phi_{i}(x_{i}) \in \mathbb{R}^{n_{i}\times d}$ is uniformly bounded in $x_i \in X_i$, i.e.,
 there exists a constant $c_{i}>0$ such that $\|\nabla\phi_{i}(x_{i})\| \leq c_{i}$ for any $x_i \in X_i$.\\
	(d)  For each agent $i\in \mathcal{N},$ $\phi_i\left(x_i\right)$ is $l_3$-Lipschitz continuous in $x_i \in X_i$.
\end{assumption}

\subsection{Graph theory}

We consider the information setting that each player $i \in \mathcal{N}$ knows the information of its private functions $f_{i},~\phi$ and the local constraint $X_i$, but have no access to the aggregate $\delta(\boldsymbol{x})$. Instead, each player is able to communicate with its neighbors over a time-varying graph $\mathcal{G}_{k}=\left\{\mathcal{N}, \mathcal{E}_{k}\right\}$. Define $W_{k}=$ $\left[\omega_{i j, k}\right]_{i, j=1}^{N}$ as the adjacency matrix of $\mathcal{G}_{k}$, where $\omega_{i j, k}>0$ if and only if $(j, i) \in \mathcal{E}_{k}$, and $\omega_{i j, k}=0$, otherwise. Denote by $N_{i, k} \triangleq\{j \in$ $\left.\mathcal{N}:(j, i) \in \mathcal{E}_{k}\right\}$ the neighboring set of player $i$ at time $k$.
We impose the following conditions on  the time-varying  communication graphs   $\mathcal{G }_k = \{ \mathcal{N}, \mathcal{E }_k\}$.

\begin{assumption}~\label{ass-graph} (a)   $W_k $ is   doubly stochastic for any $k\geq 0$;
	\\(b) There exists a constant  $0< \eta<1$  such that
	$$ \omega_{ij,k} \geq \eta  , \quad  \forall j \in  \mathcal{N}_{i,k},~~\forall i \in \mathcal{N},~\forall k \geq 0;$$
	(c)   There exists  a positive integer $B $  such that the union graph
	$\big \{ \mathcal{N }, \bigcup_{l=1}^B \mathcal{E }_{k+l} \big \}$ is strongly connected  for all $k\geq 0$.
\end{assumption}

We define a transition matrix $\Phi(k,s)  =W_kW_{k-1}\cdots W_s$ for any $k\geq s\geq 0$ with $\Phi(k,k+1)  =\mathbf{I}_N$,
and state a result    that will be used in the sequel.
\begin{lemma}\label{lemma_graph}\cite[Proposition 1]{nedic2009distributed}
	Let  Assumption \ref{ass-graph} hold. Then there exist  $\theta=(1-\eta/(4N^2))^{-2}>0$ and $\beta =(1-\eta/(4N^2))^{1/B}  $ such that for any $k\geq s \geq 0,$
	\begin{align}\label{geometric}
		\left | \left[\Phi(k,s)\right]_{ij}-1/N \right| \leq \theta \beta^{k-s},\quad \forall i,j\in \mathcal{N}.
	\end{align}
\end{lemma}

\section{Algorithm Design and Main Results }\label{sec:agg}

In this section, we design a distributed  projection-free algorithm  and provide its convergence performance.

\subsection{Distributed Projection-free Algorithm}

It is worth pointing out that the Frank-Wolfe method for the optimization problem  $\min_{x\in C} f(x) $ requires  minimizing a linear function over constraint sets \cite{scutari2012monotone}, in contrast to the projected gradient methods which require the minimization of quadratic functions over constraint sets.
Recall that the Frank-Wolfe step is given by
\begin{align}
	(F W A)\left\{\begin{array}{l}
		y_{k}=\operatorname{argmin}_{y \in C}\left\langle\nabla f\left(x_{k}\right), y\right\rangle, \\
		x_{k+1}=\left(1-\alpha_{k}\right) x_{k}+\alpha_{k} y_{k} .
	\end{array}\right.
\end{align}

Next we present the proposed  distributed Frank-Wolfe algorithm with gradient tracking (D-FWAGT). 
\begin{algorithm}[htbp]
	\caption{Distributed Projection-free Method for Aggregative Optimization}  \label{alg1}
	{\it Initialize:} Set $k=0,$ $x_{ i,0}  \in X_i$  and $ v_{ i,0}= \phi_{i}(x_{i,0})  $, and $y_{i,0} = \nabla_{v_i}g_{i}(x_{i,0},v_{i,0}) $ for  each $i \in\mathcal{N}$.
	
	{\it Iterate until convergence}\\
	{\bf Consensus.} Each  player computes an intermediate estimate by
	\begin{align}
		&\hat{v}_{i,k+1} = \sum_{j\in \mathcal{N}_{i,k}} w_{ij,k} v_{j,k} \label{alg-consensus},\\
		&\hat{y}_{i,k+1} = \sum_{j\in \mathcal{N}_{i,k}} w_{ij,k} y_{j,k} \label{alg-g}.
	\end{align}
	{\bf Strategy Update.} Each  agent $ i\in \mathcal{N}$ updates its decision  variable and its estimate of the average aggregates by
	\begin{align}
		& s_{i,k}  =\argmin_{s_{i} \in X_i} \Big( \big \langle \nabla_{1}g_{i}(x_{i,k},\hat{v}_{i,k+1})+\nabla{\phi}_{i}(x_{i,k})\hat{y}_{i,k+1},s_{i}  \big \rangle \Big) ,\notag \\
		&x_{i,k+1} = (1-\gamma_{k})x_{i,k} + \gamma_k s_{i,k} \label{alg-strategy0},\\
		& v_{i,k+1} = \hat{v}_{i,k+1}+\phi_{i}(x_{i,k+1})-\phi_{i}(x_{i,k}) ,\label{alg-average}\\
		& y_{i,k+1} = \hat{y}_{i,k+1} + \nabla_{v_{i}}g_{i}(x_{i,k+1},v_{i,k+1})-\nabla_{v_{i}}g_{i}(x_{i,k},v_{i,k}) \label{alg-gra2},
	\end{align}
	where $\gamma_{k}\in[0,1)$.
\end{algorithm}

Suppose that each agent  $i$ at stage $k $ selects a strategy  $x_{i,k} \in X_i $ as an estimate of its optimal strategy,
   and holds the estimates $v_{i,k}$ and $y_{i,k}$ for the aggregate $\delta(\boldsymbol{x})$ and the gradient term $\frac{1}{N} \sum_{i=1}^{N} \nabla_{z} g_{i}\left(x_{i}, z \right)|_{z=\delta(\boldsymbol{x})}$, respectively.
At  stage  $k+1$, player $i$   observes or receives its neighbors' information $ v_{j,k},y_{i,k} ,j\in \mathcal{N}_{i,k}$ and updates two
intermediate estimates by the consensus step \eqref{alg-consensus} and \eqref{alg-g}.
Then it computes its  gradient estimation,   and updates its strategy $x_{i,k+1} $ by a projection-free scheme  \eqref{alg-strategy0} with a Frank-Wolfe step.
Set  $y_{i,0} = \nabla_{2}g_{i}(x_{i,0},v_{i,0})  $ without loss of generality.
Finally, player $i$ updates the estimate for average aggregate $v_{i,k+1}$ with the renewed strategy  $x_{i,k+1} $  by
the dynamic average  tracking scheme  \eqref{alg-average},
and  updates the estimate of gradient  by the dynamic tracking scheme   \eqref{alg-gra2}. The procedures are summarized in  Algorithm \ref{alg1}.

Denote by  $\boldsymbol{x}_{k}:=\operatorname{col}(x_{1,k},\dots,x_{N,k})$ and similar notations for $\boldsymbol{s}_{k},~\boldsymbol{v}_{k},~\boldsymbol{y}_{k},\hat{\boldsymbol{v}}_{k },~\hat{\boldsymbol{y}}_k$. We can write the algorithm \ref{alg1} in a compact form
\begin{align}
	& \boldsymbol{s}_{k} = \argmin_{\boldsymbol{s}\in X} \langle\nabla_{1}g(\boldsymbol{x}_{k},\hat{\boldsymbol{v}}_{k+1})+\nabla{\phi}(\boldsymbol{x}_{k})\hat{\boldsymbol{y}}_{k+1},\boldsymbol{s}\rangle  \label{fw_gl} ,\\
	& \boldsymbol{x}_{k+1} = (1-\gamma_{k})\boldsymbol{x}_{k} + \gamma_{k}\boldsymbol{s}_{k} \label{fw_update},\\
	& \boldsymbol{v}_{k+1} =\hat{\boldsymbol{v}}_{k+1}   + \phi(\boldsymbol{x}_{k+1}) - \phi(\boldsymbol{x}_{k}) \label{tr_con},\\
	& \boldsymbol{y}_{k+1} =\hat{\boldsymbol{y}}_{k+1} + \nabla_{2}g(\boldsymbol{x}_{k+1},\boldsymbol{v}_{k+1}) - \nabla_{2}g(\boldsymbol{x}_{k},\boldsymbol{v}_{k})\label{tr_gra},
\end{align}
where  $\hat{\boldsymbol{v}}_{k+1}= W_{d,k}\boldsymbol{v}_{k}$, $\hat{\boldsymbol{y}}_{k+1} = W_{d,k}\boldsymbol{y}_{k}$,  $W_{d,k}=W_k\otimes I_{d}$  and $\phi(\boldsymbol{x}_{k}):=\operatorname{col}(\phi(x_{1,k}),\dots,\phi(x_{N,k})).$

We impose	the following conditions  on the step-length sequence $\left(\gamma_{k}\right)_{k \in \mathbb{N}}$.
\begin{assumption}\label{rk}
	i) (nonincreasing) $0 \leq \gamma_{k+1} \leq \gamma_{k} \leq 1$, for all $k \geq 0$;\\
	ii) (nonsummable) $\sum_{k=0}^{\infty} \gamma_{k}=\infty$;\\
	iii) (square-summable) $\sum_{k=0}^{\infty}\left(\gamma_{k}\right)^{2}<\infty$.
\end{assumption}

Denote by $\boldsymbol{x}^{\star}$ the optimal solution to the aggregative optimization problem \eqref{Nopt_agg}. We then present the main convergence result of this paper.
\begin{theorem}
	Let Algorithm 1 be applied to the    problem \eqref{Nopt_agg}, where the step size $\gamma_{k}$ satisfies  Assumption \ref{rk}. Suppose  that Assumptions \ref{ass-payoff}-\ref{ass-graph} hold. Then
	$$
	\lim_{k \rightarrow \infty} f(\boldsymbol{x}_{k}) = f(\boldsymbol{x}^{\star}).
	$$
\end{theorem}

\section{Proof of Main Result}

We first establish bounds on the consensus error of the aggregate and gradient tracking estimate  measured by
$\|\delta(\boldsymbol{x}_k) -\hat{v}_{i,k+1}\| $ and $\|\hat{\boldsymbol{y}}_{k+1} -\mathbf{1} \otimes \bar{y}_{k}  \| $, respectively, which will play an important role in the proof of Theorem 1.
Similarly to \cite[Lemma 2]{koshal2016distributed}, we have the following result.
\begin{lemma}\label{lem_ave}
	Let  Assumption \ref{ass-graph} hold. Then there exist
	\begin{align}
		& \bar{v}_{k} = \frac{1}{N}\sum_{i=1}^{N}v_{i,k} = \frac{1}{N}\sum_{i=1}^{N}\phi_{i}(x_{i,k}) , \label{average_x}\\
		& \bar{y}_{k} = \frac{1}{N}\sum_{i=1}^{N}y_{i,k} = \frac{1}{N}\sum_{i=1}^{N}\nabla_{v_{i}}g_{i}(x_{i,k},v_{i,k}). \label{average_y}
	\end{align}
\end{lemma}
\noindent {\bf Proof.} In view of \eqref{tr_con} and double-stochasticity in Assumption \ref{ass-graph}, multiplying $\frac{1}{N} \mathbf{1}_{Nd}^{\top}$ on both sides of \eqref{tr_con} can lead to that
$$
\bar{v}_{k+1}=\bar{v}_{k}+\frac{1}{N} \sum_{i=1}^{N} \phi_{i}\left(x_{i, k+1}\right)-\frac{1}{N} \sum_{i=1}^{N} \phi_{i}\left(x_{i, k}\right)
$$
which further implies that
$$
\bar{v}_{k}-\frac{1}{N} \sum_{i=1}^{N} \phi_{i}\left(x_{i, k}\right)=\bar{v}_{0}-\frac{1}{N} \sum_{i=1}^{N} \phi_{i}\left(x_{i, 0}\right) .
$$
combining the above equality and $v_{i,0} = \phi_{i}(x_{i,0})$ yields the first assertion of this lemma.

Secondly, we prove \eqref{average_y} by induction. Since the estimates are initialized as $y_{i,0}=\nabla_{v_{i}}g_{i}(x_{i,0},v_{i,0})$ for all $i\in\mathcal{N}$, we have $\bar{y}_{0} = \frac{1}{N}\sum_{i=1}^{N}\nabla_{v_{i}}g_{i}(x_{i,0},v_{i,0})$. At step $k$, we assume that $\bar{y}_{k}= \frac{1}{N}\sum_{i=1}^{N}\nabla_{v_{i}}g_{i}(x_{i,k},v_{i,k}) $. We need to  show that relation  \eqref{average_y} holds at step $k+1$:
\begin{align}
	\bar{y}_{k+1} &= \frac{1}{N}(\mathbf{1}^{\top}\otimes I_{d}) ((W_{k}\otimes I_{d})\boldsymbol{y}_{k} + \nabla_{2}g(\boldsymbol{x}_{k+1},\boldsymbol{v}_{k+1}) - \nabla_{2}g(\boldsymbol{x}_{k},\boldsymbol{v}_{k})) \notag \\
	& = \frac{1}{N}(\mathbf{1}^{\top}\otimes I_{d}) ((W_{k}\otimes I_{d})\boldsymbol{y}_{k} + \frac{1}{N}\sum_{i=1}^{N}\nabla_{v_{i}}g_{i}(x_{i,k+1},v_{i,k+1}) - \frac{1}{N}\sum_{i=1}^{N}\nabla_{v_{i}}g_{i}(x_{i,k},v_{i,k}) \notag \\
	& = \bar{y}_{k} + \frac{1}{N}\sum_{i=1}^{N}\nabla_{v_{i}}g_{i}(x_{i,k+1},v_{i,k+1}) - \frac{1}{N}\sum_{i=1}^{N}\nabla_{v_{i}}g_{i}(x_{i,k},v_{i,k}) \notag \\
	& = \frac{1}{N}\sum_{i=1}^{N}\nabla_{v_{i}}g_{i}(x_{i,k+1},v_{i,k+1}), \notag
\end{align}
where the first equality follows by the update rule of $y_{i}'s$ in \eqref{alg-average} of Algorithm \ref{alg1}, the second follows from definition of $\bar{y}_{k}$, i.e. $\bar{y}_{k} = \frac{1}{N}(\mathbf{1}^{\top}\otimes I_d)\boldsymbol{y}_{k}$, the third follows by Assumption \ref{ass-graph}. While the last equality follows from the induction step $k$. \hfill $\blacksquare$

Then, we  establish a  bound on the consensus error $\|\delta(\mathbf{x}_k) -\hat{v}_{i,k+1}\|$ of the aggregate. This proof is similar to \cite{koshal2016distributed}, we state the proof of this proposition based on algorithm (\ref{alg1}) to  ensure the integrity of the work.
For each $i\in \mathcal{N},$ we  define
\begin{align}
	M_i & \triangleq   \max_{x_i \in X_i} \|x_i\|,~ M_H\triangleq \sum_{j=1}^N M_i, \label{def-bdst} \\
	& {\rho}_i:=\max_{\theta,\theta^{\prime}\in X_i}\|\theta-\theta^{\prime}\|, \rho \triangleq \max_{i\in \mathcal{N}}\rho_{i} .\label{x_bound}
\end{align}

\begin{proposition} \label{prop_agg} Consider Algorithm \ref{alg1}. Let Assumptions \ref{ass-payoff}, and \ref{ass-graph}  hold.  Then
	\begin{equation}\label{agg-bd0}
		\|\delta(\boldsymbol{x}_k) -\hat{v}_{i,k+1}\|\leq \theta M_H \beta^{k}    +\theta N \rho l_3  \sum_{s=1}^k \beta^{k-s}   \gamma_{s-1},
	\end{equation}
	where the constants $\theta, \beta$ are defined in \eqref{geometric}.
\end{proposition}
\noindent {\bf Proof.}   From \eqref{average_x} it follows  that
\begin{align}\label{equi}  \sum_{i=1}^N v_{i,k} =\sum_{i=1}^N  \phi_i (x_{i,k}) = N\delta(\boldsymbol{x}_k),\quad \forall k\geq 0.
\end{align}

Akin to   \cite[Eqn. (16)]{koshal2016distributed}, we give an upper bound on
$\left \|{ \delta(\boldsymbol{x}_k)} -\hat{v}_{i,k+1} \right \|.$
By combining \eqref{alg-average}  with \eqref{alg-consensus}, we have
\begin{align*}  v_{i,k+1}  &= \sum_{j=1}^N[ \Phi (k,0)]_{ij} v_{j,0 }
	+ \phi_{i}(x_{i,k+1})- \phi_{i}(x_{i,k})
	+\sum_{s=1}^{k }  \sum_{j=1}^N [ \Phi (k,s)]_{ij}( \phi_{j}(x_{j,s})- \phi_{i}(x_{j,s-1}) ) .
\end{align*}
Then  by \eqref{alg-average}, we have
$$ \hat{v}_{i,k+1 }
=\sum_{j=1}^N[ \Phi (k,0)]_{ij} v_{j,0 } +\sum_{s=1}^{k }  \sum_{j=1}^N [ \Phi (k,s)]_{ij}( \phi_{j}(x_{j,s})- \phi_{j}(x_{j,s-1}) )  . $$
By  using \eqref{equi}, we have that
\begin{align*}
	{ \delta(\boldsymbol{x}_k)}= {\sum_{j=1}^N v_{j,0}\over N}+ \sum_{s=1}^k \sum_{j=1}^N {1\over N}(  \phi_{j}(x_{j,s})- \phi_{j}(x_{j,s-1}) )  .\end{align*}
Therefore,  we obtain the following bound.
\begin{align*}
	\left \|{ \delta(\boldsymbol{x}_k)} -\hat{v}_{i,k+1} \right \| &\leq \sum_{j=1}^N \left| {1\over N}-[\Phi(k,0)]_{ij}\right| \|v_{j,0}\|  \notag\\
	&+\sum_{s=1}^k \sum_{j=1}^N
	\left| {1\over N}-[\Phi(k,s)]_{ij}\right|  \big\|   \phi_{j}(x_{j,s})- \phi_{j}(x_{j,s-1})   \big \|.
\end{align*}
Then by using \eqref{geometric}, and  $v_{i,0}= x_{i,0} $,  we obtain that
\begin{equation}\label{bd-consensus1}
	\begin{split}
		\left \|{ \delta(\boldsymbol{x}_k)} -\hat{v}_{i,k+1} \right  \| \leq  \theta \beta^{k}\sum_{j=1}^N \| x_{j,0} \|   +   \theta \sum_{s=1}^k \beta^{k-s}   \sum_{j=1}^N   \big\|  \phi_{j}(x_{j,s})- \phi_{j}(x_{j,s-1})  \big \| .
\end{split} \end{equation}

  Note that for any $s \ge 1$ and each $i\in \mathcal{N}$,
\begin{align}\label{phi}
	\|  \phi_{i}(x_{i,s})- \phi_{i}(x_{i,s-1})\| \leq l_3\| x_{i,s}- x_{i,s-1}\| = l_3 \gamma_{s-1}\| s_{i,s-1}- x_{i,s-1}\| \leq  l_3 \gamma_{s-1}\rho,
\end{align}
where the first inequality follows by Assumption \ref{ass-gradient}(d) (i.e., the $l_3$-Lipschitz continuous of $\phi_{i}$), and the second equality follows by \eqref{alg-strategy0}, and the last inequality follows by the constraint set is convex and bounded as \eqref{x_bound}. By combining \eqref{phi}, \eqref{bd-consensus1} and \eqref{def-bdst}, we  prove  \eqref{agg-bd0}.
\hfill $\blacksquare$

Next, we establish a bound on the consensus error  $\|\hat{\boldsymbol{y}}_{k+1} -\mathbf{1} \otimes \bar{y}_{k}  \| $ of the gradient tracking step, where $\bar{y}_{k}$ is defined by  Lemma \ref{lem_ave}.

\begin{proposition} \label{Pro_gra}
	Consider Algorithm \ref{alg1}. Suppose that  Assumption \ref{ass-payoff}-\ref{ass-graph} hold. Define $C_k \triangleq  \max_{i\in \mathcal{N}}\left \|{ \delta(\boldsymbol{x}_k)} -\hat{v}_{i,k+1} \right\|$. Then the following hold for all $k\in \mathbb{N}$:
	\begin{align}\label{average_gra}
		\|\hat{\boldsymbol{y}}_{k+1} -\mathbf{1} \otimes \bar{y}_{k}  \| \leq & \theta \beta^{k}\|\boldsymbol{y}_{0} \| + \sum_{s=1}^{k}\theta\beta^{k-s}(l_2\bar{\rho } + l_2l_3\sqrt{N}\rho) \gamma_{s-1} \notag \\
		& +  \sum_{s=1}^{k} l_2 l_3 (\sqrt{N}+1)\rho\theta\beta^{k-s} \gamma_{s-2}+\sum_{s=1}^{k}l_2\sqrt{N}\theta\beta^{k-s}(C_{s-1} +C_{s-2}),
	\end{align}
\end{proposition}
where $\theta, \beta$ are defined in \eqref{geometric},  $\bar{\rho }$ is the diameter of convex set $X$ as in Assumption \ref{ass-payoff}.

\noindent{\bf{Proof:}}  
By telescoping \eqref{tr_gra}, we obtain
\begin{align}
	\boldsymbol{y}_{k+1} &= W_{d,k}(W_{d,k-1}\boldsymbol{y}_{k-1} + \nabla_2 g(\boldsymbol{x}_{k},\boldsymbol{v}_{k}) - \nabla_2 g(\boldsymbol{x}_{k-1},\boldsymbol{v}_{k-1})) + \nabla_2 g(\boldsymbol{x}_{k+1},\boldsymbol{v}_{k+1}) - \nabla_2 g(\boldsymbol{x}_{k},\boldsymbol{v}_{k}) \notag \\
	& =\dots \notag \\
	& = (\Phi(k,0) \otimes I_d) \boldsymbol{y}_{0} + \sum_{s=1}^{k}(\Phi(k,s)\otimes I_d)(\nabla_2 g(\boldsymbol{x}_{s},\boldsymbol{v}_{s}) - \nabla_2 g(\boldsymbol{x}_{s-1},\boldsymbol{v}_{s-1})) \notag\\
	& + \nabla_2 g(\boldsymbol{x}_{k+1},\boldsymbol{v}_{k+1}) - \nabla_2 g(\boldsymbol{x}_{k},\boldsymbol{v}_{k}). \label{y_inter}
\end{align}
By arranging \eqref{tr_gra}, we can write $	 W_{d,k} \boldsymbol{y}_{k} =\boldsymbol{y}_{k+1} - \nabla_2 g(\boldsymbol{x}_{k+1},\boldsymbol{v}_{k+1}) + \nabla_2 g(\boldsymbol{x}_{k},\boldsymbol{v}_{k})$. Then by exploiting the equivalence in \eqref{y_inter}, we have
\begin{align}\label{y_estimate}
	W_{d,k} \boldsymbol{y}_{k} = (\Phi(k,0) \otimes I_d) \boldsymbol{y}_{0} + \sum_{s=1}^{k}(\Phi(k,s) \otimes I_d)(\nabla_2 g(\boldsymbol{x}_{s},\boldsymbol{v}_{s}) - \nabla_2 g(\boldsymbol{x}_{s-1},\boldsymbol{v}_{s-1}))
\end{align}

By equation \eqref{average_y}, we have $\bar{y}_{s} = \frac{1}{N}(\mathbf{1}^{\top}\otimes I_d)(\nabla_2 g(\boldsymbol{x}_{s},\boldsymbol{v}_{s})),\forall s\geq 0,$ which leads to:
\begin{align}\label{y_average}
	\bar{y}_{k} = \frac{1}{N}(\mathbf{1}^{\top}\otimes I_d)\boldsymbol{y}_0 + \sum_{s=1}^{N}\frac{1}{N}(\mathbf{1}^{\top}\otimes I_d)(\nabla_2 g(\boldsymbol{x}_{s},\boldsymbol{v}_{s})-\nabla_2 g(\boldsymbol{x}_{s-1},\boldsymbol{v}_{s-1})).
\end{align}
From \eqref{y_estimate} and \eqref{y_average}, we have the following:
\begin{align}\label{bd-y1}
	&\|W_{d,k}\boldsymbol{y}_{k} -\mathbf{1}\otimes \bar{y}_{k} \| \notag \\ &=\|(\Phi(k,0)-\frac{1}{N}\mathbf{1}\mathbf{1}^{\top})\boldsymbol{y}_{0}+\sum_{s=1}^{k}(\Phi(k,s)-\frac{1}{N}\mathbf{1}\mathbf{1}^{\top})(\nabla_2 g(\boldsymbol{x}_{s},\boldsymbol{v}_{s}) - \nabla_2 g(\boldsymbol{x}_{s-1},\boldsymbol{v}_{s-1}))\| \notag \\
	& \leq \|\Phi(k,0)-\frac{1}{N}\mathbf{1}\mathbf{1}^{\top} \| \|\boldsymbol{y}_{0} \|+ \sum_{s=1}^{k}\|\Phi(k,s)-\frac{1}{N}\mathbf{1}\mathbf{1}^{\top}\| \|\nabla_2 g(\boldsymbol{x}_{s},\boldsymbol{v}_{s}) - \nabla_2 g(\boldsymbol{x}_{s-1},\boldsymbol{v}_{s-1})\|  \\
	& \leq \theta \beta^{k}\|\boldsymbol{y}_{0} \| + \sum_{s=1}^{k}\theta\beta^{k-s}\|\nabla_2 g(\boldsymbol{x}_{s},\boldsymbol{v}_{s}) - \nabla_2 g(\boldsymbol{x}_{s-1},\boldsymbol{v}_{s-1})\| ,\notag
\end{align}
where the first inequality follows from the Kronecker and  Cauchy-Schuarz inequality, and the second inequality follows by Lemma \ref{lemma_graph}.

Next, we find an upper bound for $\|\nabla_2 g(\boldsymbol{x}_{s},\boldsymbol{v}_{s}) - \nabla_2g(\boldsymbol{x}_{s-1},\boldsymbol{v}_{s-1})\|$.  Note by \eqref{phi} that
\begin{align*}
&  \|\phi(\boldsymbol{x}_{s})-\phi(\boldsymbol{x}_{s-1})\| =\sqrt{\sum_{i=1}^N 	\|  \phi_{i}(x_{i,s})- \phi_{i}(x_{i,s-1})\|^2} \leq   \sqrt{N} l_3\rho \gamma_{s-1},
\\& \|\delta(\boldsymbol{x}_{s-1})-\delta(\boldsymbol{x}_{s-2})\|= \frac{1}{N}\left \|\sum_{i=1}^{N}
(\phi_{i}(x_{i,s-1})- \sum_{i=1}^{N}\phi_{i}(x_{i,s-2})) \right \| \leq   l_3 \rho\gamma_{s-2}.
\end{align*}
Then by recalling  that  $\bar{\rho}:=\max_{\theta,\theta^{\prime}\in X}\|\theta-\theta^{\prime}\|^{2}_{2}$, and using \eqref{fw_update} and \eqref{tr_con}, we have:
\begin{align}\label{bd-g}
	&\|\nabla_2 g(\boldsymbol{x}_{s},\boldsymbol{v}_{s}) - \nabla_2g(\boldsymbol{x}_{s-1},\boldsymbol{v}_{s-1})\| \leq l_2\|\boldsymbol{x}_{s}-\boldsymbol{x}_{s-1}\|+l_2\|\boldsymbol{v}_{s}-\boldsymbol{v}_{s-1} \| \notag\\
	&= l_2 \| \gamma_{s-1}\boldsymbol{s}_{s-1}-\gamma_{s-1}\boldsymbol{x}_{s-1}\| \notag\\
	&+ l_2\|\hat{\boldsymbol{v}}_{s}  + \phi(\boldsymbol{x}_{s}) - \phi(\boldsymbol{x}_{s-1})-\mathbf{1}\otimes \delta(\boldsymbol{x}_{s-1})+\mathbf{1}\otimes \delta(\boldsymbol{x}_{s-1})-\mathbf{1}\otimes \delta(\boldsymbol{x}_{s-2})\notag \\
	&+\mathbf{1}\otimes \delta(\boldsymbol{x}_{s-2})-(\hat{\boldsymbol{v}}_{s-1}  + \phi(\boldsymbol{x}_{s-1}) - \phi(\boldsymbol{x}_{s-2})) \|\notag \\
	& \leq l_2 \bar{\rho } \gamma_{s-1} + l_2\|\hat{\boldsymbol{v}}_{s} - \mathbf{1}\otimes\delta(\boldsymbol{x}_{s-1})\|+ l_2\|\phi(\boldsymbol{x}_{s})-\phi(\boldsymbol{x}_{s-1})\| + l_2\|\delta(\boldsymbol{x}_{s-1})-\delta(\boldsymbol{x}_{s-2})\|\notag\\
	& + l_2\| \mathbf{1}\otimes\delta(\boldsymbol{x}_{s-2})-\hat{\boldsymbol{v}}_{s-1}\| +l_2\|\phi(\boldsymbol{x}_{s-1})-\phi(\boldsymbol{x}_{s-2})\|\notag\\
	& \leq l_2\bar{\rho}\gamma_{s-1} + l_2\sqrt{N} C_{s-1} +l_2l_3\sqrt{N}\rho\gamma_{s-1} +l_2l_3\rho\gamma_{s-2} +l_2\sqrt{N} C_{s-2}+ l_2l_3\sqrt{N}\rho\gamma_{s-2} \notag \\
	& \leq l_2\bar{\rho}\gamma_{s-1}  +l_2l_3\sqrt{N}\rho\gamma_{s-1} +  l_2l_3(\sqrt{N}+1)\rho\gamma_{s-2} + l_2\sqrt{N} C_{s-1} +l_2\sqrt{N} C_{s-2},
\end{align}
where the first inequality follows from  Assumption \ref{ass-gradient}(b) (i.e., the $l_2$-lipschitz of continuous of $\nabla_2 g(\boldsymbol{x},\boldsymbol{z})$), the last inequality  but one holds by definition   $C_k = \max_{i\in \mathcal{N}}\left \|{ \delta(\boldsymbol{x}_k)} -\hat{v}_{i,k+1} \right\|$ and  \eqref{phi}. Finally, by substituting \eqref{bd-g} into \eqref{bd-y1} we  derive \eqref{average_gra}. \hfill $\blacksquare$

Now, we state a convergence property of Algorithm \ref{alg1}.
\begin{proposition} \label{prop_FW} Consider Algorithm \ref{alg1}. Let Assumptions \ref{ass-payoff}-\ref{ass-graph}  hold.  Define  $h_{k}=f(\boldsymbol{x}_{k})-f(\boldsymbol{x}^{\star})$ and   $c^{\prime}=\max_{i\in \mathcal{N}}c_{i}$. Then
	\begin{equation}\label{convergence-pre}
		h_{k+1}\leq (1-\gamma_{k})h_k + \frac{L}{2}\gamma_{k}^{2}\bar{\rho}^{2}    +{\varepsilon_{1,k}}\gamma_{k} +{\varepsilon_{2,k}}\gamma_{k},
	\end{equation}
	where
	\begin{align}
		\varepsilon_{1,k} &= \bar{\rho}l_1\sqrt{N} (\theta \beta^{k-1}M_H  +   \theta N \rho \sum_{s=1}^k \beta^{k-s} r_{s-1}) ,\label{def-varek1} \\
		\varepsilon_{2,k} &= \rho c^{\prime}\big(\theta \beta^{k}\|\boldsymbol{y}_{0} \| + \sum_{s=1}^{k}\theta\beta^{k-s}(l_2\bar{\rho } + l_2l_3\sqrt{N}\rho) \gamma_{s-1}  +  \sum_{s=1}^{k}l_2l_3(\sqrt{N}+1) \rho\theta\beta^{k-s} \gamma_{s-2}\notag \\
		& +\sum_{s=1}^{k}l_2\sqrt{N}\theta\beta^{k-s}(C_{s-1} +C_{s-2}) +\sqrt{N}l_2 C_{k-1}+2\sqrt{N}l_2l_3 \rho \gamma_{k-1}\big) . \label{def-varek2}
	\end{align}
	
\end{proposition}

\noindent{\bf{Proof:}} Note by \eqref{fw_update} that
 $\boldsymbol{x}_{k+1}-\boldsymbol{x}_{k}  =  \gamma_{k}(\boldsymbol{s}_{k}-\boldsymbol{x}_{k} )$.
Then from the L-smoothness of $f$ and the boundedness of $X$, we have
\begin{align} \label{pro_main}
	f(\boldsymbol{x}_{k+1}) &\leq f(\boldsymbol{x}_{k}) +\langle \nabla f(\boldsymbol{x}_{k}),\boldsymbol{x}_{k+1}-\boldsymbol{x}_{k}  \rangle +\frac{L}{2}\|\boldsymbol{x}_{k+1}-\boldsymbol{x}_{k}\|^{2} \notag \\
	&\leq f(\boldsymbol{x}_{k}) +\gamma_{k}\langle \nabla f(\boldsymbol{x}_{k}),\boldsymbol{s}_{k}-\boldsymbol{x}_{k}  \rangle +\frac{L}{2}\gamma_{k}^{2}\bar{\rho}^{2},
\end{align}
where $\bar{\rho}:=\max_{\theta,\theta^{\prime}\in X}\|\theta-\theta^{\prime}\|^{2}_{2}$.

Note that
\begin{align}
	\nabla f(\boldsymbol{x}_{k}) = \nabla_{1}g(\boldsymbol{x}_{k},\boldsymbol{\delta}(\boldsymbol{x}_{k}))+\nabla\phi(\boldsymbol{x}_{k})\mathbf{1}\otimes\frac{1}{N}\sum_{i=1}^{N}\nabla_2g_{i}(x_{i,k},\delta(\boldsymbol{x}_{k}))
\end{align}
where $\boldsymbol{\delta}(\boldsymbol{x}_{k}) \triangleq \mathbf{1}\otimes \delta(\boldsymbol{x}_{k}) $. Thus, we see that for any  $\boldsymbol{s}\in X$,
\begin{align*}
	\langle 	\nabla f(\boldsymbol{x}_{k}),\boldsymbol{s}  \rangle &= \langle \nabla_{1}g(\boldsymbol{x}_{k},\hat{\boldsymbol{v}}_{k+1})+ \nabla\phi(\boldsymbol{x})\hat{\boldsymbol{y}}_{k+1}, \boldsymbol{s} \rangle  +\langle \nabla_{1}g(\boldsymbol{x}_{k},\boldsymbol{\delta}(\boldsymbol{x})) - \nabla_{1}g(\boldsymbol{x}_{k},\hat{\boldsymbol{v}}_{k+1}), \boldsymbol{s} \rangle \notag \\
	&+ \left \langle \nabla\phi(\boldsymbol{x}_{k})\mathbf{1}\otimes\frac{1}{N}\sum_{i=1}^{N}\nabla_2g_{i}(x_{i},\delta(\boldsymbol{x})) - \nabla\phi(\boldsymbol{x}_{k})\hat{\boldsymbol{y}}_{k+1}, \boldsymbol{s}   \right \rangle.
\end{align*}
Thus, we have
\begin{align}\label{err_k}
	\langle 	\nabla f(\boldsymbol{x}_{k}),\boldsymbol{s}_{k}-\boldsymbol{x}^{\star} \rangle &= \langle \nabla_{1}g(\boldsymbol{x}_{k},\hat{\boldsymbol{v}}_{k+1})+ \nabla\phi(\boldsymbol{x}_{k})\hat{\boldsymbol{y}}_{k+1}, \boldsymbol{s}_{k}-\boldsymbol{x}^{\star}\rangle \notag\\
	& +\langle \nabla_{1}g(\boldsymbol{x},\boldsymbol{\delta}(\boldsymbol{x})) - \nabla_{1}g(\boldsymbol{x}_{k},\hat{\boldsymbol{v}}_{k+1}), \boldsymbol{s}_{k}-\boldsymbol{x}^{\star}\rangle \notag \\
	&+ \left \langle \nabla\phi(\boldsymbol{x}_{k})\mathbf{1}\otimes\frac{1}{N}\sum_{i=1}^{N}\nabla_2g_{i}(x_{i,k},\delta(\boldsymbol{x}_{k})) - \nabla\phi(\boldsymbol{x}_{k})\hat{\boldsymbol{y}}_{k+1}, \boldsymbol{s}_{k}-\boldsymbol{x}^{\star}  \right \rangle \notag \\
	&\leq \langle \nabla_{1}g(\boldsymbol{x}_{k},\boldsymbol{\delta}(\boldsymbol{x}_{k})) - \nabla_{1}g(\boldsymbol{x}_{k},\hat{\boldsymbol{v}}_{k+1}), \boldsymbol{s}_{k}-\boldsymbol{x}^{\star}\rangle \notag \\
	&+ \left \langle \nabla\phi(\boldsymbol{x}_{k})\mathbf{1}\otimes\frac{1}{N}\sum_{i=1}^{N}\nabla_2g_{i}(x_{i,k},\delta(\boldsymbol{x}_{k})) - \nabla\phi(\boldsymbol{x}_{k})\hat{\boldsymbol{y}}_{k+1}, \boldsymbol{s}_{k}-\boldsymbol{x}^{\star}  \right \rangle
\\& \leq \bar{\rho}\|\nabla_{1} g(\boldsymbol{x}_{k},\boldsymbol{\delta}(\boldsymbol{x}_{k}))-\nabla_{1} g(\boldsymbol{x}_{k},\hat{\boldsymbol{v}}_{k+1})\| \notag \\
	&+\bar{\rho} \left \|\nabla \phi(\boldsymbol{x}_k)\mathbf{1}_{N}\otimes\frac{1}{N}\sum_{i=1}^{N}\nabla_{2}g_{i}(x_{i,k},\delta(\boldsymbol{x}_k))-\nabla \phi(\boldsymbol{x}_k)\hat{\boldsymbol{y}}_{k+1} \right \| \notag,
\end{align}
where  the first inequality hold since  $\langle \nabla_{1}g(\boldsymbol{x}_{k},\hat{\boldsymbol{v}}_{k+1})+ \nabla\phi(\boldsymbol{x}_{k})\hat{\boldsymbol{y}}_{k+1}, \boldsymbol{s}_{k}-\boldsymbol{x}^{\star}\rangle \leq 0$ by \eqref{fw_gl},  and the last inequality has   utilized Assumption \ref{ass-payoff}(a).
 Then by adding $\langle \nabla f(\boldsymbol{x}_{k}),\boldsymbol{x}^{\star}-\boldsymbol{x}_{k} \rangle$ to both side of \eqref{err_k} and using the triangle inequality, we have
\begin{align}\label{h_k}
	\langle \nabla f(\boldsymbol{x}_{k}),\boldsymbol{s}_{k}-\boldsymbol{x}_{k} \rangle
	& \leq \langle \nabla f(\boldsymbol{x}_{k}),\boldsymbol{x}^{\star}-\boldsymbol{x}_{k} \rangle +\bar{\rho}\|\nabla_{1} g(\boldsymbol{x}_{k},\boldsymbol{\delta}(\boldsymbol{x}_{k}))-\nabla_{1} g(\boldsymbol{x}_{k},\hat{\boldsymbol{v}}_{k+1})\| \notag \\
	&+\bar{\rho} \| \nabla \phi(\boldsymbol{x}_k) \| \left \| \mathbf{1}_{N}\otimes\frac{1}{N}\sum_{i=1}^{N}\nabla_{2}g_{i}(x_{i,k},\delta(\boldsymbol{x}_k))- \hat{\boldsymbol{y}}_{k+1} \right \| \notag \\
	&\leq \langle \nabla f(\boldsymbol{x}_{k}),\boldsymbol{x}^{\star}-\boldsymbol{x}_{k} \rangle +\bar{\rho}l_{1}\|\boldsymbol{\delta}(\boldsymbol{x}_{k})-\hat{\boldsymbol{v}}_{k+1}\| \notag \\ &+\bar{\rho}c^{\prime}\left \|\mathbf{1}\otimes\frac{1}{N}\sum_{i=1}^{N}\nabla_{2}g_{i}(x_{i,k},\delta(\boldsymbol{x}_k))
-\hat{\boldsymbol{y}}_{k+1}\right \|,
\end{align}
where the last inequality has applied   Assumption \ref{ass-gradient}(a),  Assumption \ref{ass-gradient}(b), and  $c^{\prime}=\max_{i\in \mathcal{N}}c_{i}$.

Next, we prove the bound of {$\|\mathbf{1}\otimes\frac{1}{N}\sum_{i=1}^{N}\nabla_{2}g_{i}(x_{i},\delta(\boldsymbol{x}))-\hat{\boldsymbol{y}}_{k+1}\| $. By adding and subtracting  $\frac{1}{N}\sum_{i=1}^{N}\nabla_{2}g_{i}(x_{i,k},v_{i,k})$, we have
	\begin{align*}
		& \left \|\hat{\boldsymbol{y}}_{k+1} - \mathbf{1} \otimes \frac{1}{N}\sum_{i=1}^{N}\nabla_{2}g_{i}(x_{i,k},\delta(\boldsymbol{x}_{k})) \right \|
\notag \\
		& =\left  \|\hat{\boldsymbol{y}}_{k+1} -\mathbf{1} \otimes \frac{1}{N}\sum_{i=1}^{N}\nabla_{2}g_{i}(x_{i,k},v_{i,k}) +\mathbf{1}\otimes\frac{1}{N}\sum_{i=1}^{N}\nabla_{2}g_{i}(x_{i,k},v_{i,k})- \mathbf{1}\otimes\frac{1}{N}\sum_{i=1}^{N}\nabla_{2}g_{i}(x_{i,k},\delta(\boldsymbol{x}_{k})) \right \| \notag  \\
		& {\eqref{average_y} \atop =} \left  \|\hat{\boldsymbol{y}}_{k+1} -\mathbf{1} \otimes \bar{y}_{k} \right \| + \left \|\mathbf{1}\otimes \frac{1}{N}\sum_{i=1}^{N}\nabla_{2}g_{i}(x_{i,k},v_{i,k})- \mathbf{1}\otimes\frac{1}{N}\sum_{i=1}^{N}\nabla_{2}g_{i}(x_{i,k},\delta(\boldsymbol{x}_k)) \right \| \notag \\
		&\leq \|\hat{\boldsymbol{y}}_{k+1} -\mathbf{1} \otimes \bar{y}_{k} \| + \left \|\mathbf{1}\otimes\frac{1}{N}\sum_{i=1}^{N}\nabla_{2}g_{i}(x_{i,k},v_{i,k})-
\mathbf{1}\otimes\frac{1}{N}\sum_{i=1}^{N}\nabla_{2}g_{i}(x_{i,k},\delta(\boldsymbol{x}_{k-1}))\right\| \notag \\
		& +\left \|\mathbf{1}\otimes\frac{1}{N}\sum_{i=1}^{N}\nabla_{2}g_{i}(x_{i,k},\delta(\boldsymbol{x}_{k-1}))- \mathbf{1}\otimes\frac{1}{N}\sum_{i=1}^{N}\nabla_{2}g_{i}(x_{i,k},\delta(\boldsymbol{x}_{k})) \right \| \notag .
	\end{align*}
	Then by  using the triangle inequality  and Assumption \ref{ass-gradient}, we obtain that
	\begin{align}
		&\left  \|\hat{\boldsymbol{y}}_{k+1} - \mathbf{1} \otimes \frac{1}{N}\sum_{i=1}^{N}\nabla_{2}g_{i}(x_{i,k},\delta(\boldsymbol{x}_{k})) \right \| \notag\\
		& \leq \|\hat{\boldsymbol{y}}_{k+1} -\mathbf{1} \otimes \bar{y}_{k} \|  + \frac{1}{\sqrt{N}}\sum_{i=1}^{N}\|\nabla_{2}g_{i}(x_{i,k},v_{i,k})-\nabla_{2}g_{i}(x_{i,k},\delta(\boldsymbol{x}_{k-1}))\|  \notag \\
		&+\frac{1}{\sqrt{N} }\sum_{i=1}^{N} \|\nabla_{2}g_{i}(x_{i,k},\delta(\boldsymbol{x}_{k}))-\nabla_{2}g_{i}(x_{i,k},\delta(\boldsymbol{x}_{k-1}))\|\notag\\
		&\leq  \|\hat{\boldsymbol{y}}_{k+1} -\mathbf{1} \otimes \bar{y}_{k} \|+\frac{1}{\sqrt{N}}\sum_{i=1}^{N}l_2\|\delta(\boldsymbol{x}_{k-1}) -v_{i,k} \|  + \sqrt{N}l_2\|\delta(\boldsymbol{x}_{k})-\delta(\boldsymbol{x}_{k-1})\| \notag
		\\&{\eqref{alg-average} \atop =} \|\hat{\boldsymbol{y}}_{k+1} -\mathbf{1} \otimes \bar{y}_{k} \|+\sqrt{N}l_2\|\delta(\boldsymbol{x}_{k})-\delta(\boldsymbol{x}_{k-1})\|\notag\\
		&+\frac{1}{\sqrt{N}}\sum_{i=1}^{N}l_2\|\delta(\boldsymbol{x}_{k-1}) -\hat{v}_{i,k} + \phi_i(x_{i,k-1}) -\phi_i(x_{i,k})\|  \notag \\
		&\leq \|\hat{\boldsymbol{y}}_{k+1} -\mathbf{1} \otimes \bar{y}_{k} \|+\sqrt{N}l_2 C_{k-1} +\sqrt{N}l_2 l_3\rho\gamma_{k-1} +\sqrt{N}l_2l_3\rho \gamma_{k-1},
	\end{align}
	where the last inequality holds by  the definition $C_k = \max_{i\in \mathcal{N}}\left \|{ \delta(\boldsymbol{x}_k)} -\hat{v}_{i,k+1} \right\|$  and \eqref{phi}. Then by applying Proposition \ref{Pro_gra} and recalling the definition of
	$\varepsilon_{2,k}$ in \eqref{def-varek2}, we have
	\begin{align}\label{grad_error}
		\bar{\rho}c^{\prime}\|\hat{\boldsymbol{y}}_{k+1} - \mathbf{1} \otimes \frac{1}{N}\sum_{i=1}^{N}\nabla_{2}g_{i}(x_{i,k},\delta(\boldsymbol{x}_k))\| \leq \varepsilon_{2,k}.
	\end{align}

	Note by Proposition \ref{prop_agg} that
	\begin{align}\label{N_agg_error}
	\bar{\rho}l_1	\|\boldsymbol{\delta}(\boldsymbol{x}_{k})-\hat{\boldsymbol{v}}_{k+1}\| \leq \bar{\rho}l_1 \sqrt{N}(\theta \beta^{k-1}M_H  +   \theta N \rho \sum_{s=1}^k \beta^{k-s} r_{s-1}){\eqref{def-varek1}\atop =}\varepsilon_{1,k}.
	\end{align}
	Then by substituting \eqref{grad_error} and \eqref{N_agg_error} into  \eqref{h_k}, we derive
	\begin{align*}
		\langle \nabla f(\boldsymbol{x}_{k}),\boldsymbol{s}_{k}-\boldsymbol{x}_{k} \rangle \leq \langle \nabla f(\boldsymbol{x}_{k}),\boldsymbol{x}^{\star}-\boldsymbol{x}_{k} \rangle +\varepsilon_{1,k}+ \varepsilon_{2,k} .\end{align*}
	Therefore, by subtracting $f(\boldsymbol{x}^{\star})$ from both sides  in the inequality \eqref{pro_main}, and using the above inequality, we obtain that
	\begin{align} \label{conv_ineq}
		h_{k+1} \leq h_{k} +\gamma_{k}\langle\boldsymbol{x}^{\star}-\boldsymbol{x}_{k}, \nabla f(\boldsymbol{x}_{k}) \rangle + \frac{L}{2}\gamma_{k}^{2}\bar{\rho}^{2}  +{\varepsilon_{1,k}}\gamma_{k} +{\varepsilon_{2,k}}\gamma_{k}.
	\end{align}
	Since $f$ is convex, we observe $
		\langle\boldsymbol{x}^{\star}-\boldsymbol{x}_{k}, \nabla f(\boldsymbol{x}_{k}) \rangle \leq -h_k.$
This incorporating with \eqref{conv_ineq} proves the proposition. \hfill $\blacksquare$
	
	Next, we give the following convergence result of a recursive linear inequality.
	\begin{proposition}\label{sequence} \cite[Lemma 3]{polyak1987introduction}  Let the nonnegative sequence $\left\{u_{k}\right\}$ be generated by $u_{k+1} \leq q_{k} u_{k} + \nu_{k}$, where $0 \leq q_{k} \leq 1, \nu_{k} \geq 0$.
		Suppose that
		\begin{align}\label{sequence_polyak}
\sum_{k=0}^{\infty}\left(1-q_{k}\right)=\infty, \quad \nu_{k} /\left(1-q_{k}\right) \rightarrow 0.
		\end{align}
		Then $\lim _{k \rightarrow \infty} u_{k} \leq 0$. In particular, if $u_{k} \geq 0$, then, $u_{k} \rightarrow 0$
	\end{proposition}
	
	\begin{lemma}\label{infty} [\cite{sundhar2010distributed}, Lemma 3]
		Let $\left(\gamma_{k}\right)_{k \in \mathbb{N}}$ be a scalar sequence.\\
		a) If $\lim _{k \rightarrow \infty} \gamma_{k}=\gamma$ and $0<\tau<1$, then $\lim _{k \rightarrow \infty}$ $\sum_{\ell=0}^{k} \tau^{k-\ell} \gamma_{\ell}=\gamma /(1-\tau)$;\\
		b) If $\gamma_{k} \geq 0$ for all $k, \sum_{k=0}^{\infty} \gamma_{k}<\infty$ and $0<\tau<1$, then $\sum_{k=0}^{\infty} \sum_{\ell=0}^{k} \tau^{k-\ell} \gamma_{\ell}<\infty .$
	\end{lemma}

	Based on Propositions \ref{prop_agg}, \ref{Pro_gra}, \ref{prop_FW} and \ref{sequence}, we proceed to prove that
	the iterates generated by the proposed algorithm \ref{alg1} converge to the optimal point.

	\noindent{\bf{Proof of Theorem 1:}} Set  $q_k = 1- \gamma_{k}$, and $\nu_k = \frac{L}{2}\gamma_{k}^{2}\bar{\rho}^{2} +\varepsilon_{1,k}\gamma_{k} + \varepsilon_{2,k}\gamma_{k}$. We will apply Proposition 4 to  \eqref{convergence-pre} to prove  Theorem 1.
	
	Firstly, we prove that \[\lim_{k\rightarrow \infty} \varepsilon_{1,k} + \varepsilon_{2,k} = 0.\]

	Since $\lim_{k\rightarrow \infty}\gamma_{k} = 0$ by Assumption  \ref{rk} and  $0 < \beta < 1 $
by Lemma \ref{lemma_graph},  we drive $\lim_{k\rightarrow \infty}\sum_{s=1}^k \beta^{k-s}\gamma_{s-1} =0 $ and $\lim_{k\rightarrow \infty}\sum_{s=1}^k \beta^{k-s}\gamma_{s-2} =0 $ by Lemma \ref{infty}(a). Thus,
	\begin{align}
		\lim_{k\rightarrow \infty}\varepsilon_{1,k} = \lim_{k\rightarrow \infty}~\bar{\rho}l_1\sqrt{N} (\theta \beta^{k-1}M_H  +   \theta N \rho \sum_{s=1}^k \beta^{k-s} r^{s-1})  = 0.
	\end{align}
	Then by definition $C_k = \max_{i\in \mathcal{N}}\left \|{ \delta(\boldsymbol{x}_k)} -\hat{v}_{i,k+1} \right\|$  and using Proposition \ref{prop_agg}, we have $\lim_{k\rightarrow \infty} C_k = 0$.
	Similarly, by  using Lemma \ref{infty}(a),  we obtain that
	\[ \lim_{k\rightarrow \infty}\sum_{s=1}^k \beta^{k-s} C_{s-1} =0 ~{\rm and} ~\lim_{k\rightarrow \infty}\sum_{s=1}^k \beta^{k-s} C_{s-2} = 0.\] Thus,
	\begin{align}
		\lim_{k\rightarrow \infty}\varepsilon_{2,k} &= \lim_{k\rightarrow \infty} \rho c^{\prime}\big(\theta \beta^{k}\|\boldsymbol{y}_{0} \| + \sum_{s=1}^{k}\theta\beta^{k-s}(l_2\bar{\rho } + l_2l_3\sqrt{N}\rho) \gamma_{s-1}  +  \sum_{s=1}^{k}l_2l_3(\sqrt{N}+1) \rho\theta\beta^{k-s} \gamma_{s-2}\notag \\
		& +\sum_{s=1}^{k}l_2\sqrt{N}\theta\beta^{k-s}(C_{s-1} +C_{s-2}) +\sqrt{N}l_2 C_{k-1}+
(\sqrt{N}+1) l_2l_3 \rho \gamma_{k-1}\big) = 0.
	\end{align}
	
From Assumption \ref{rk} it follows that $\sum_{0}^{k}(1 - q_{k}) =\sum_{0}^{k}\gamma_{k} = \infty$ and $\lim_{k\rightarrow \infty} \frac{L\bar{\rho}^2}{2}\gamma_{k} = 0$.
	In summary, we have proved the condition \eqref{sequence_polyak} required by Proposition \ref{sequence}. Then by applying Proposition \ref{sequence} to \eqref{convergence-pre}, we obtain $\lim_{k\rightarrow \infty} f(\boldsymbol{x}_{k})-f(\boldsymbol{x}^{\star})=0$. \hfill $\blacksquare$
	
	\section{Numerical Simulation}\label{sec:Num}
	
	In this section,  we demonstrate the proposed algorithm by
	solving an example with $N=5$ agents for problem (\ref{Nopt_agg}).  Agent $i$'s local cost function is
	\begin{align}\label{Simulation}
	  f_{i}\left(x_{i}, \delta(x)\right)=k_{i}\left(x_{i}-\chi_{i}\right)^{2}+P(\delta(x)) x_{i}
	\end{align}
    where $k_{i}$ is constant and $\chi_{i}$ is the fixed entities for $i=1, \cdots, N$, and $P(\sigma(x))=$ $a N \sigma(x)+p_{0}$ with $\sigma(x)=\frac{1}{N} \sum_{i \in \mathcal{N}} x_{i}$.
	In simulation, the constraint set is set as an $l_{1}$ norm ball constraint $\Omega_i=\left\{x_i \mid\|x_i\|_{1} \leq R_i \right\}$. Then, $s_{i,k}$ in D-FWAGT admits a closed form solution
	$$
	s_{i,k}=\operatorname{argmin}_{s \in \Omega_i}\left\langle s, d_{i,k}\right\rangle=R_i \cdot(-\operatorname{sgn}\left[d_{i,k}\right])
	$$
	with  $d_{i,k} \triangleq \nabla_{1}g_{i}(x_{i,k},\hat{v}_{i,k+1})+\nabla{\phi}_{i}(x_{i,k})\hat{y}_{i,k+1} $ as in algorithm \ref{alg1}.
    Let the feasible sets is $\Omega_i=\left\{x_i \in \mathbb{R}^{n}: \|x_1\|_{1} \leq 5,\right.$ $\left.\|x_2\|_{1} \leq 7,\right.$ $\left.\|x_3\|_{1} \leq 9, \right.$ $\left.\|x_4\|_{1} \leq 3,\right.$ $\left.\|x_5\|_{1} \leq 6 \right\}$.  And we have $col(\chi_{i})_{i=1,\dots,5} = [3,5,6,1,2]^{\top}\otimes \mathbf{1}_n$, $a=0.04$, and  $p_{0} =5\times\mathbf{1}_n.$
	
	 We set an undirected time-varying graph as the communication network. The graph at each iteration is randomly drawn from a set of three graphs, whose union graph is connected. Set the adjacency matrix $W=\left[w_{i j}\right]$, where $w_{i j}=\frac{1}{\max \left\{\left|\mathcal{N}_{i }\right|,\left|\mathcal{N}_{j}\right|\right\}}$ for any $i \neq j$ with $(i, j) \in \mathcal{E}, w_{i i}=1-\sum_{j \neq i} w_{i j}$, and $w_{i j}=0$.
	
	\begin{figure}[htbp]
		\centering
		\begin{subfigure}[$\gamma_k = 1/k, \gamma_k =1/\sqrt{k}$]
			{
				\begin{minipage}{0.3\linewidth}
					\includegraphics[width=0.9\textwidth, height=0.25\textheight]{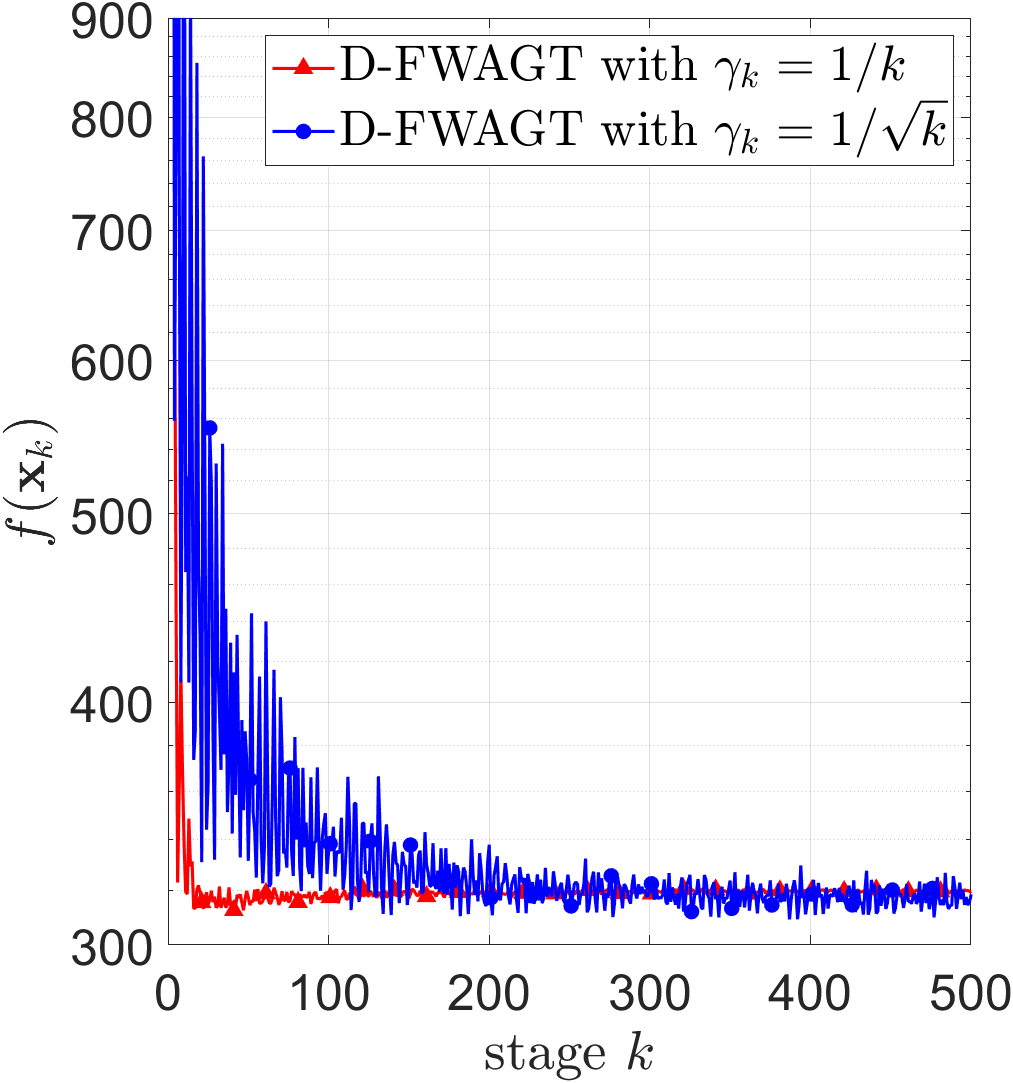}
				\end{minipage}
			}
		\end{subfigure}
		\centering
		\begin{subfigure}[$\gamma_k = 1/k, \gamma_k =1/k^2$]
			{
				\begin{minipage}{0.3\linewidth}
					\includegraphics[width=0.9\textwidth, height=0.25\textheight]{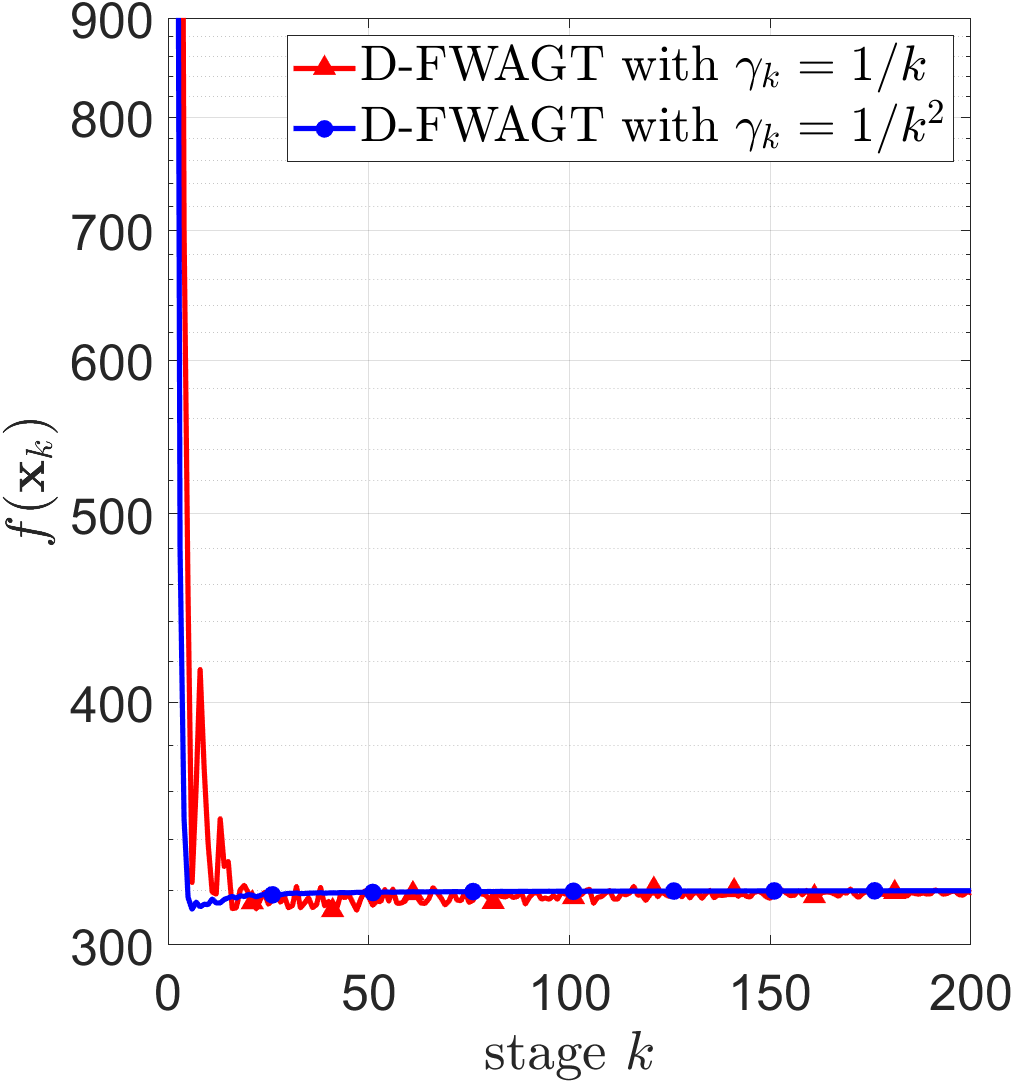}
				\end{minipage}
			}
		\end{subfigure}
		\caption{Evolution of $f(\boldsymbol{x})'s$ versus iteration $s$.}
		\label{fig1}
	\end{figure}

    Fig.\ref{fig1} displays the convergence of the proposed algorithm, and it compares the  effect of the step size of the Frank Wolfe type algorithms on the convergence of the algorithm, by considering $\gamma_k = 1/\sqrt{k},$ $\gamma_k = 1/{k},$ and $\gamma_k = 1/{k}^2.$  Although the convergence rate of the algorithm is almost as fast for the three different decreasing step cases, too large decreasing step size could cause the result of the algorithm to be too far from the optimal solution.

\begin{figure}[H]
	\centering
	\begin{subfigure}[$n=32$]
		{
		\begin{minipage}{0.3\linewidth}
		\includegraphics[width=0.9\textwidth, height=0.25\textheight]{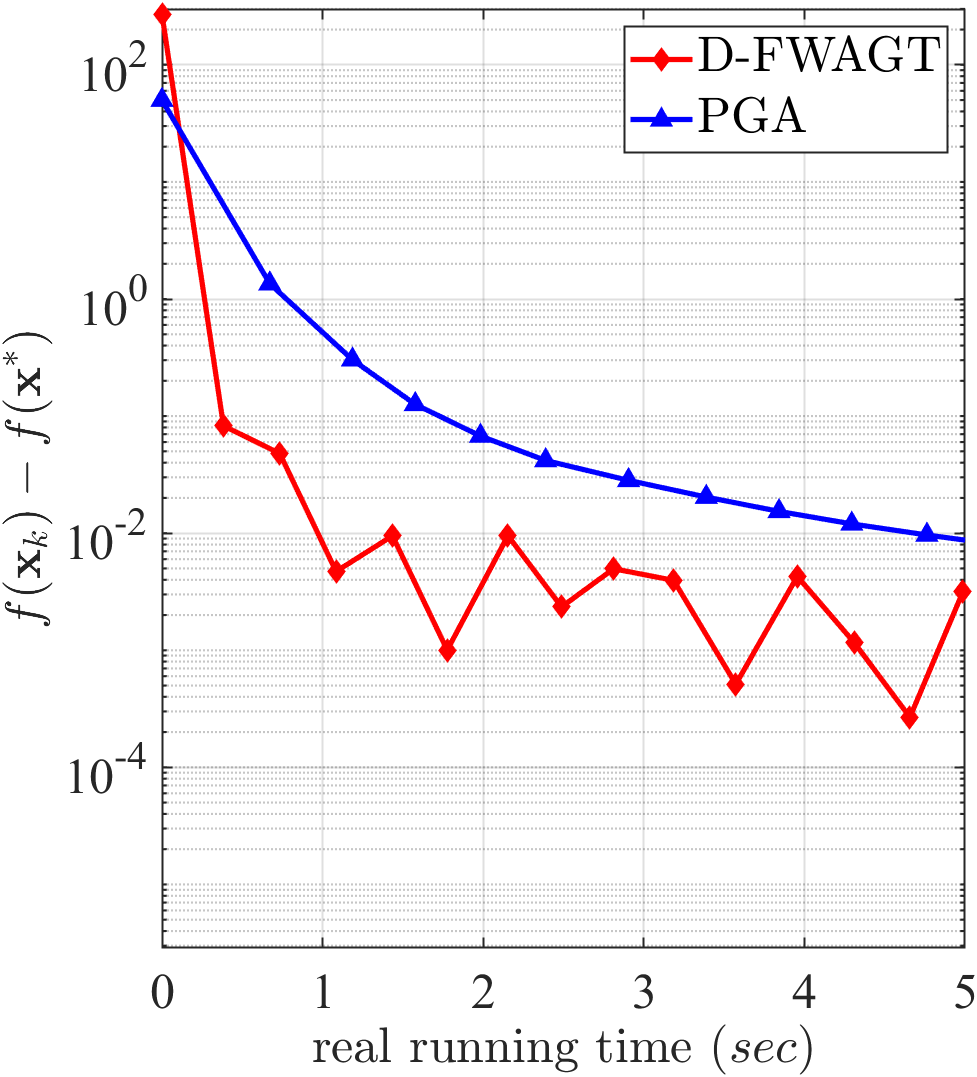}
		\end{minipage}
	}
	\end{subfigure}
    \centering
	\begin{subfigure}[$n=64$]
		{
			\begin{minipage}{0.3\linewidth}
				\includegraphics[width=0.9\textwidth, height=0.25\textheight]{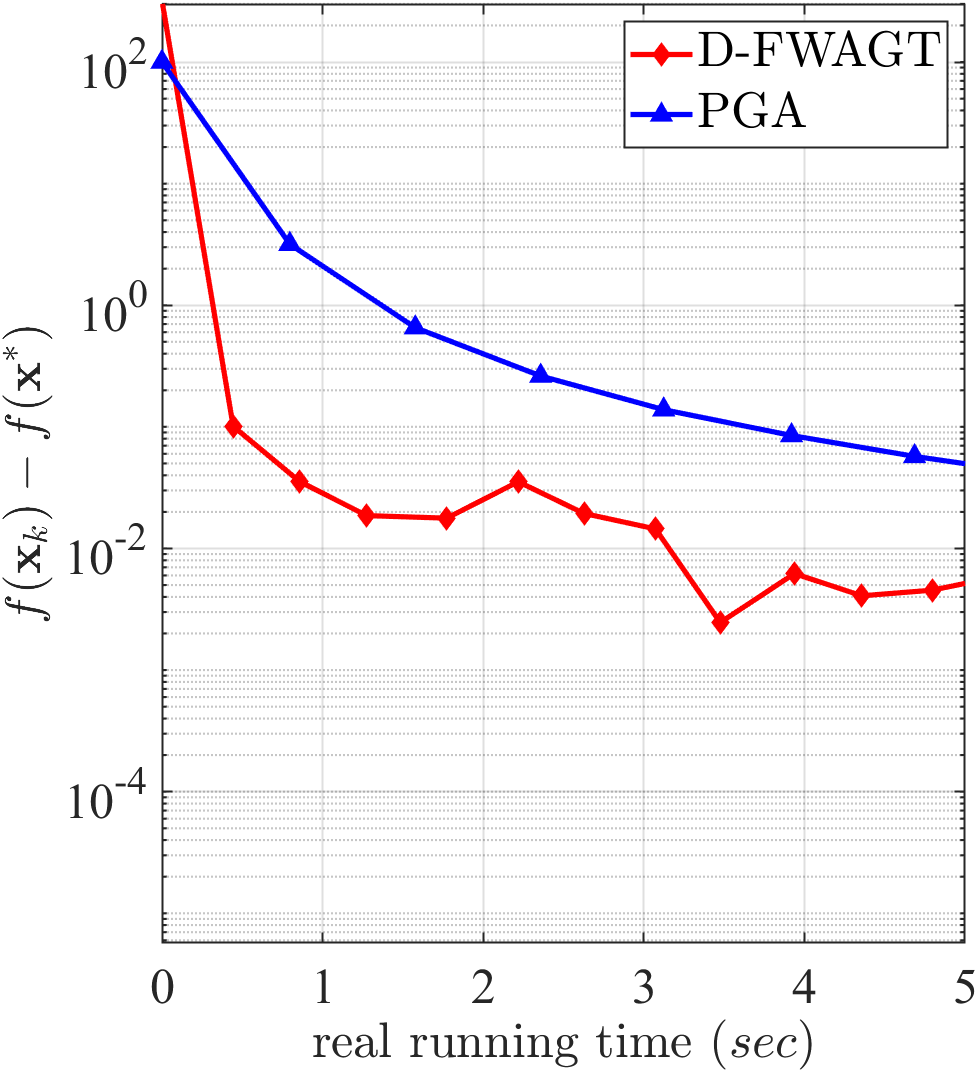}
			\end{minipage}
		}
	\end{subfigure}
\centering
	\begin{subfigure}[$n=128$]
		{
			\begin{minipage}{0.3\linewidth}
				\includegraphics[width=0.9\textwidth, height=0.25\textheight]{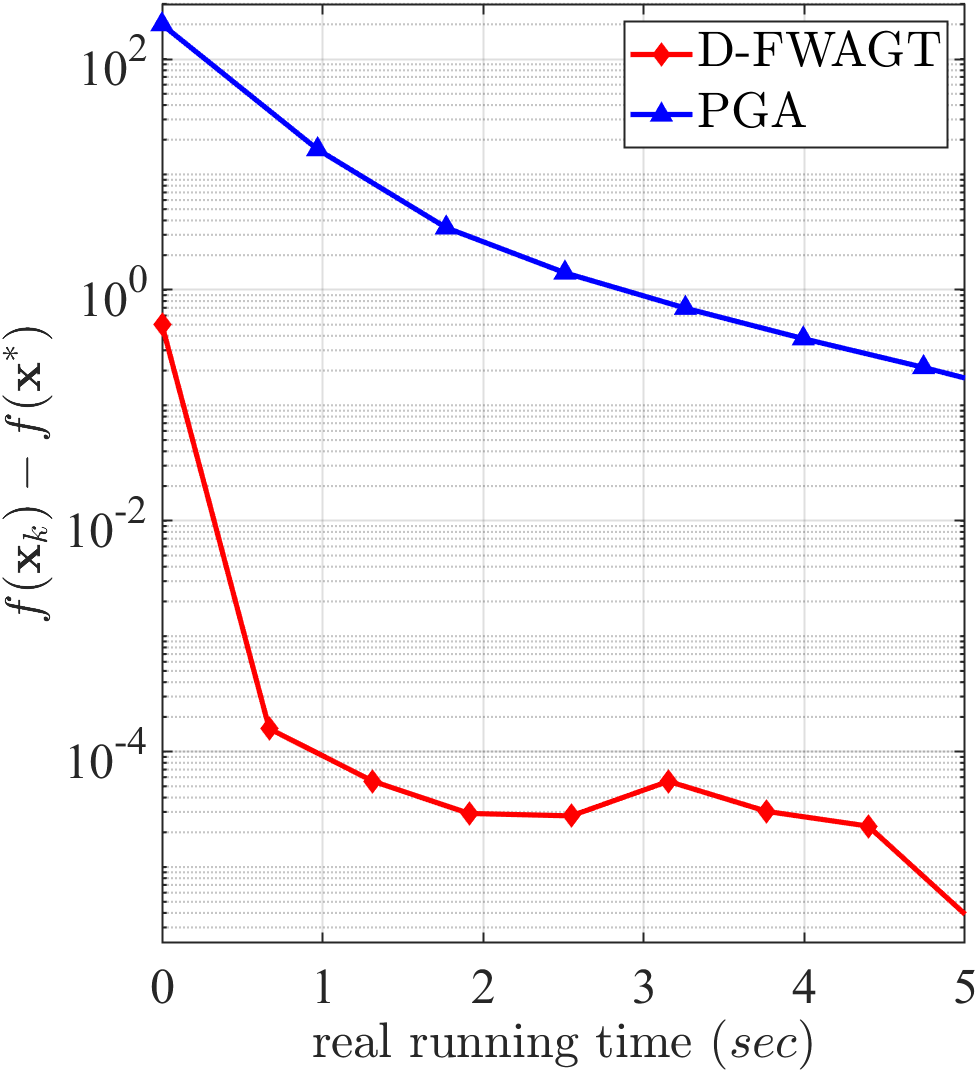}
			\end{minipage}
		}
	\end{subfigure}
\caption{Optimal solution errors with different dimensions n = 32, 64, 128.}
\label{fig2}
\end{figure}

To demonstrate the properties of our algorithm. We compare our algorithm \ref{alg1} with projection based algorithm ($x_{i,k+1} =\Pi_{x_i}(x_{i,k}-\alpha \nabla_{1}g_{i}(x_{i,k},\hat{v}_{i,k+1})+\nabla{\phi}_{i}(x_{i,k})\hat{y}_{i,k+1})$), and select the dimensions of decision variables as the power of two. In Fig.\ref{fig2}, the x-axis is for the real running time (CPU time) in seconds, while the y-axis is for the optimal solution errors in each algorithm. We learn from Fig.\ref{fig2} that as the dimension increases, the actual running time (CPU time) of the projection-based algorithm is significantly longer than that of the projection-free algorithm. The reason is that searching for poles on the boundary of a high-dimensional constraint set (solving a linear program) is faster than computing the projection of a high-dimensional constraint set (solving a quadratic program).

	\begin{table}[H]
		\centering
		\setlength{\tabcolsep}{5mm}
	    \begin{tabular}{c|c|c|c|c|c}
		\hline \hline
		\text { dimensions } & {n}=16 & {n}=32 & {n}=64 & {n}=128 & {n}=256  \\
		\hline \text { D-FWAGT }({msec}) & 0.069 & 0.077 & 0.092 & 0.126 & 0.132  \\
		\hline \text { PGA (msec) } & 78.5 &95.6  & 164.6   & 242.2 & 366.3   \\
		\hline \hline
	    \end{tabular}
    \caption{THE AVERAGE REAL RUNNING TIME OF SOLVING ONE-STAGE SUBPROBLEMS.}
    \label{tab1}
	\end{table}
	
	In addition, in Tab. \ref{tab1}, we list the average actual running time for solving the single-stage subproblem, i.e., linear or quadratic programs. When the dimensional is low, the difference in time required to solve the linear program and the quadratic program on such constraint sets may not be too great. However, as the dimension explodes, solving the quadratic program becomes difficult in this case, but the time to solve the linear program hardly varies much. That is consistent with the advantages of the projection free approaches for large-scale problems.

	\section{Conclusions }\label{sec:con}
	
	This paper proposes  a distributed projection free gradient method for aggregative optimization problem based on Frank-Wolfe method, and shows that the proposed method can achieve the convergence for the case  of the cost function is convex. In addition, empirical results demonstrate that our method indeed brings speed-ups. It is of interest to explore the faster convergence rate projection-free algorithm  for distributed aggregative optimization, and analysis the Frank-Wolfe method to the other classes of network optimizations in distributed and stochastic settings.

\bibliographystyle{elsarticle-num}
\bibliography{references}           

\end{document}